\documentclass[a4,10pt,leqno]{amsart}
\usepackage{graphicx}
\usepackage{enumerate}
\newcommand{\red}[1]{#1}

\newtheorem{theorem}{Theorem}[section]
\newtheorem{lemma}[theorem]{Lemma}
\newtheorem{proposition}[theorem]{Proposition}
\newtheorem{cororally}[theorem]{Cororally}
\newtheorem{definition}[theorem]{Definition}
\newtheorem{remark}[theorem]{Remark}
\newtheorem{example}[theorem]{Example}
\newtheorem{acutually}[theorem]{Fact}
\newtheorem{exercise}{Exercise}[section]

\newenvironment{lem}{\begin{lemma}}{\end{lemma}}

\newenvironment{df}{\begin{definition}\rm}{\end{definition}}
\newenvironment{rem}{\begin{remark}\rm}{\end{remark}}
\newenvironment{ex}{\begin{example}\rm}{\end{example}}
\newenvironment{fact}{\begin{acutually}}{\end{acutually}}

\makeatletter

\@addtoreset{equation}{section}
\makeatother

  \newcommand{\fulltoday}{\ifcase\month\or
    January\or February\or March\or April\or May\or June\or
    July\or August\or September\or October\or November\or December\fi
    \space\number\day,\space \number\year}

\newcommand{\R}{\mathbb{R}}

\newcommand{\Sb}{\mathbb{S}}

\newcommand{\ve}{\boldsymbol{e}}
\newcommand{\vn}{\boldsymbol{n}}
\newcommand{\vb}{\boldsymbol{b}}

\title[Singularities of developable M\"obius strips]
{Singularities of the asymptotic completion of developable M\"obius strips}

\author[K.~Naokawa]{Kosuke~Naokawa}

\date{\fulltoday}

\address{
      Department of Mathematics,
      Graduate School of Science,
      Osaka University,
      Toyonaka, Osaka 560-0043, Japan
}

\email{k-naokawa@cr.math.sci.osaka-u.ac.jp}


\begin{document}

\begin{abstract}
   We prove that the asymptotic completion of a developable M\"obius strip in Euclidean three-space
   must have at least one singular point other than cuspidal edge singularities.
   Moreover, if the strip contains a closed geodesic,
   then the number of such singular points is at least three.
   These lower bounds are both sharp. 
\end{abstract}

\maketitle

\section{Introduction}

%Research of flat M\"obius strips in Euclidean space start 

%Flat M\"obius strips in Euclidean 3-space were made up 

Let $U$ be an open domain in Euclidean two-space $\R^2$ and
$f: U \longrightarrow \R^3$ a $C^\infty$ map.
A  point $p\in U$ is called a \textit{singular point} of $f$
if the Jacobi matrix of $f$ is of rank less than $2$ at $p$.
It is well-known that complete and flat
 (i.e.\ zero Gaussian curvature) 
surfaces immersed in $\R^3$
are cylindrical.
This fact implies that the `asymptotic completion'
(see Definition \ref{def:a-completion})
of a developable 
M\"obius strip (i.e.\ a flat ruled M\"obius strip)
must have singular points.
%where a singular point of
%a $C^\infty$ map $f$ of a domain in $\R^2$ into $\R^3$
%is a point where $f$ is not an immersion.
Since the most generic singular points appeared on
developable
surfaces are cuspidal edge singularities (cf.\ \cite{Izumiya,Izumiya2}), 
we are interesting how often singular points other than 
cuspidal edge singularities appear on the asymptotic completion 
of a developable M\"obius strip. 
(Izumiya-Takeuchi \cite{Izumiya2} is a nice reference for
singularities of ruled surfaces or developable surfaces.)

Recently, global properties of flat surfaces 
with singularities in $\R^3$ 
were investigated in Murata-Umehara \cite{Murata}.
They defined `completeness' for flat fronts (cf.\ \cite[Definition 0.2]{Murata})
and proved that a complete flat front
with embedded ends has at least four singular 
points other than cuspidal edge singularities
if the front has singular points.
However, we cannot apply this result, since 
complete flat fronts are all orientable (cf.\ \cite[Theorem A]{Murata}).
Therefore, it is interesting to determine lower bounds on
the number of non-cuspidal-edge singular points
on developable M\"obius strips.
We show the following:

\medskip
\noindent
\textbf{Proposition.}\ \it
%\begin{prop} %\label{prop:mainA}
The asymptotic completion of a developable M\"obius strip has
at least one singular point other than cuspidal edge singularities.
\rm %\end{prop}
\medskip

In fact, there are many developable M\"obius strips.
Chicone-Kalton \cite{Chicone} constructed
a developable M\"obius strip on each generic closed regular curve
in $\R^3$.
The topological types of M\"obius strips are determined by
the isotopy types of their generating curves
and M\"obius twisting numbers.
R{\o}gen \cite{Rogen} showed that
there exists a developable M\"obius strip of an arbitrarily 
given topological type.

A developable M\"obius strip which contains a closed geodesic 
is called a \textit{rectifying M\"obius strip}.
Roughly speaking, a rectifying strip can be 
constructed from an isometric deformation of a rectangular domain
on a plane (cf.\ \cite[Proposition 2.14]{Kurono}% or Sabitov{S}
).
We also show the following assertion.

\medskip
\noindent
\textbf{Theorem.}\ \it
%\begin{thm} \label{thm:mainB}
The asymptotic completion of a rectifying M\"obius strip
has at least three singular points other than cuspidal edge singularities.
%\end{thm}
\rm 

\medskip
The first explicit construction of a rectifying M\"obius strip in $\R^3$
was given by Wunderlich \cite{Wunderlich}.
Recently, Kurono-Umehara \cite{Kurono} proved that
there exists a rectifying  M\"obius strip %(resp.\ principal M\"obius strip)
which is isotopic to any given M\"obius strip.
See Sabitov \cite{Sabitov}
for other references and the history.

\section{Singularities of developable M\"obius strips}
First, we define several terminologies.
Let $\gamma = \gamma(s) : \R \longrightarrow \R^3$
be a $C^\infty$ map.
The map $\gamma(s)$ is called 
\textit{$l$-periodic} if $\gamma(s+l)=\gamma(s)$ for $s \in \R$.
A $C^\infty$ map $\gamma(s)$ is called \textit{regular} if
$\gamma'(s):=d\gamma(s)/ds$ does not vanish on $\R$.
We fix such a periodic regular curve $\gamma$.
An $\R^3$-valued vector field $\xi$ along $\gamma$
is called \textit{$l$-odd-periodic} if it satisfies $\xi(s+l)=-\xi(s)$ for $s \in \R$.
We also fix such a $C^\infty$ $l$-odd-periodic vector field $\xi$.
Then, a $C^\infty$ immersion
\begin{equation}\label{eq:F}
      F(s, u)=\gamma(s) + u \xi(s) \qquad (s\in\R,\ |u|<\epsilon)
\end{equation}
is called a \textit{ruled M\"obius strip} if 
$\gamma'(s)$ and $\xi(s)$ are linearly independent for each $s\in \R$,
where $\epsilon>0$ is taken to be sufficiently small.
In this situation, $\gamma$ is called the \textit{generating curve} of $F$
and $\xi$ is called the \textit{ruling vector field} of $F$.

\begin{df}\label{def:a-completion}
   Let $F(s,u)$ be a ruled M\"obius strip as in \eqref{eq:F}.
   Then a $C^\infty$ map
$$
\tilde{F}(s,u)=\gamma(s) + u \xi(s) \qquad (s,u\in \R)
$$
is called the \textit{asymptotic completion} (or \textit{a-completion}) of $F$.
\end{df}

%Next, we define the following equivalence relation to investigate types of singular points.

\begin{figure}[htbp]
   \begin{minipage}{0.49\hsize}
      \begin{center}
         \includegraphics[width=0.35\hsize]{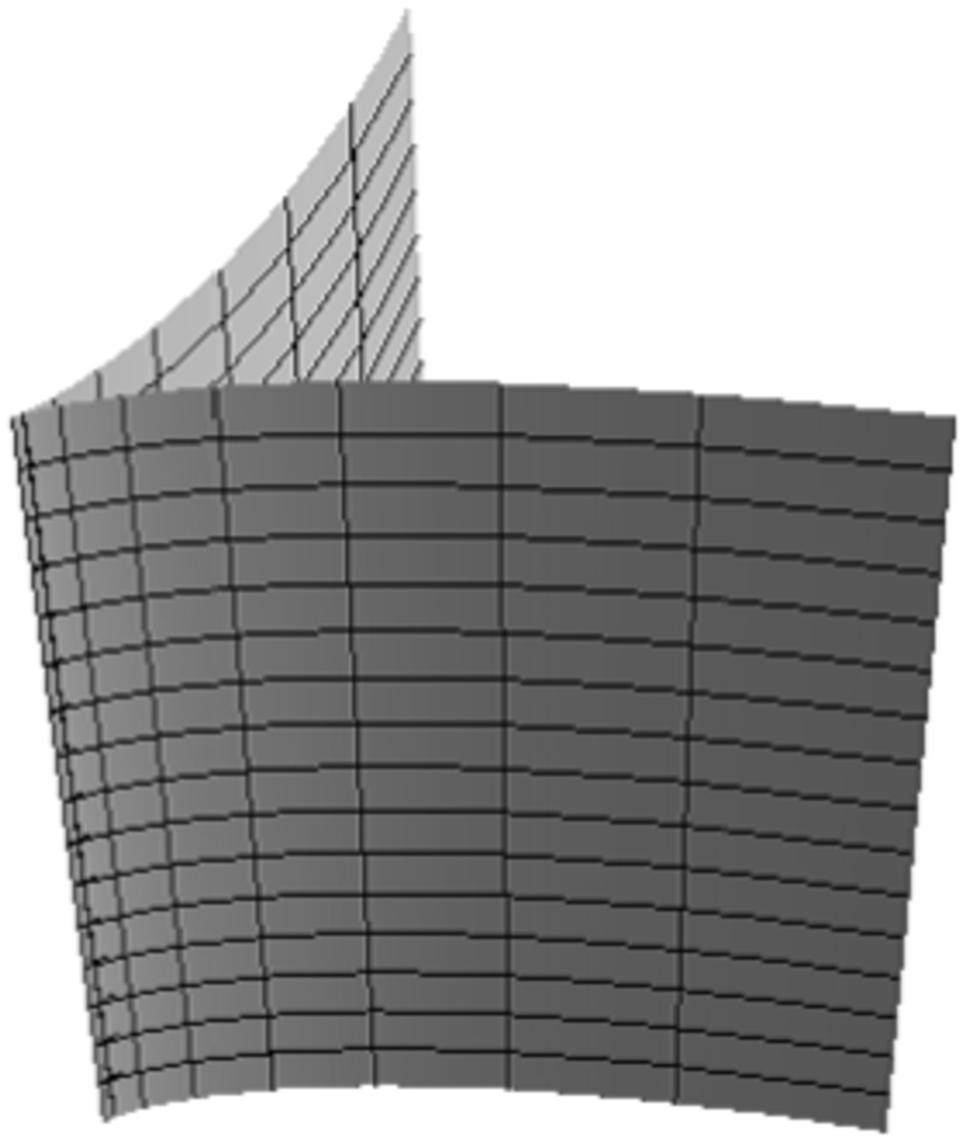}
      \end{center}
      \caption{A cuspidal edge}
      \label{fig:cuspidal_edge}
   \end{minipage}
   \begin{minipage}{0.49\hsize}
      \begin{center}
         \includegraphics[width=0.6\hsize]{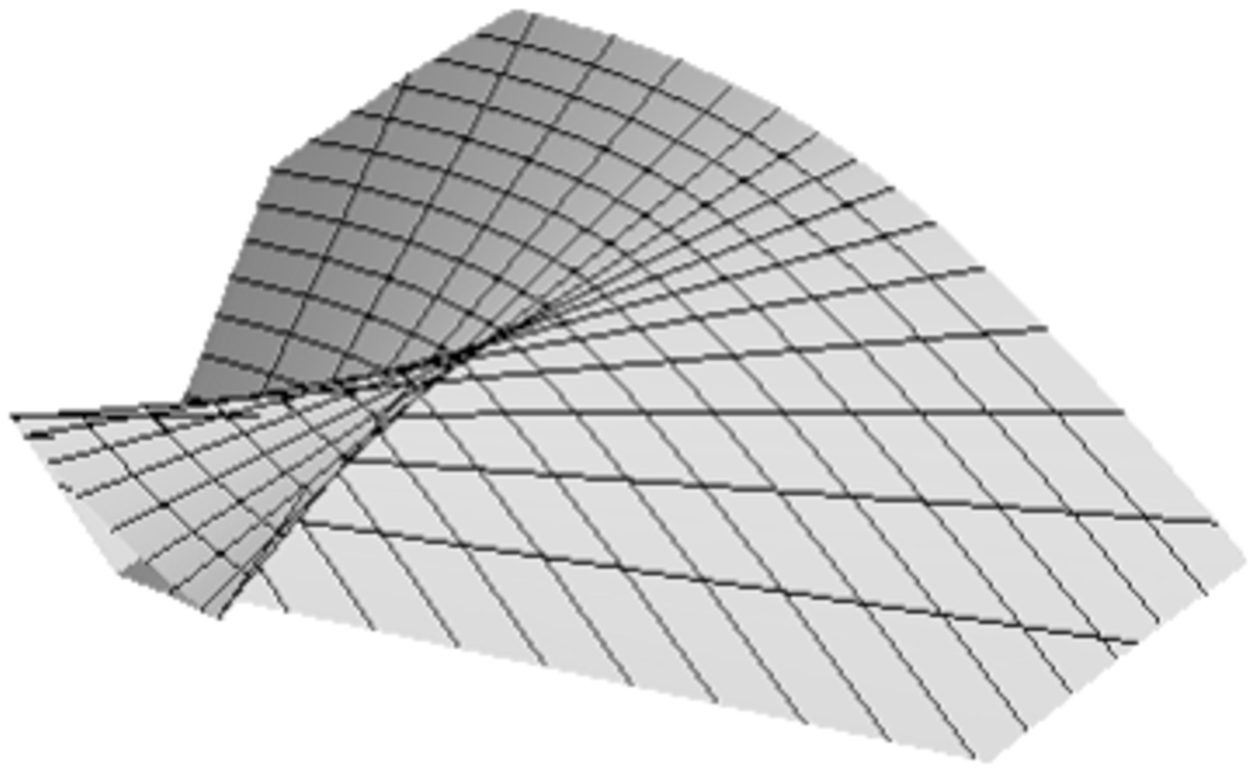}
      \end{center}
      \caption{A swallowtail}
      \label{fig:swallowtail}
   \end{minipage}
\end{figure}

\begin{df}
Let $U_i \subset \R^2$ ($i=1,2$)
be two open neighborhoods of points $p_i \in \R^2$ and
$f_i = f_i(u,v): U_i \longrightarrow \R^3$ two $C^\infty$ maps.
Then $f_1$ is said to be \textit{right-left equivalent} to $f_2$
if there exist two diffeomorphisms $\varphi : \R^2 \longrightarrow \R^2$
and $\Phi : \R^3 \longrightarrow \R^3$ such that 
$\varphi(p_1)=p_2$, $\Phi\circ f_1(p_1)=f_2(p_2)$ and $\Phi \circ f_1 = f_2 \circ \varphi$ on $U_1$.
\end{df}
We set
$$
f_C(u,v):=\begin{pmatrix}
2u^3\\  -3u^2\\  v
\end{pmatrix},
\qquad
f_S(u,v):=\begin{pmatrix}
3u^4+u^2v \\ -4u^3-2u v \\ v
\end{pmatrix}
$$
(See Figures \ref{fig:cuspidal_edge} and \ref{fig:swallowtail}, respectively).
A $C^\infty$ map which is right-left equivalent to the map germ
$f_C$ (resp.\ $f_S$)
is called a \textit{cuspidal edge} (resp.\ a \textit{swallowtail}).

%and
%the map $f_S(u,v):=(3u^4+u^2v, -4u^3-2u v, v)$ is called a \textit{cuspidal edge}
%
%
%A map germ which is right-left equivalent in a neighborhood of the origin to
%the map $f_C(u,v):=(2u^3, -3u^2, v)$
%that is the Cartesian product of a $3/2$-cusp and real interval and has the singular point at the origin
%is said to the \textit{cuspidal edge}.
%A map germ which is right-left equivalent in a neighborhood of the origin to
%the map $f_S(u,v):=(3u^4+u^2v, -4u^3-2u v, v)$
%is said to the \textit{swallowtail}.

%\begin{figure}[htbp]
%   \begin{minipage}{0.49\hsize}
%         \includegraphics[width=0.35\hsize]{cuspidal_edge.eps}
%      \caption{A cuspidal edge}
%      \label{fig:cuspidal_edge}
%   \end{minipage}
%   \begin{minipage}{0.49\hsize}
%%      \begin{center}
%         \includegraphics[width=0.6\hsize]{swallowtail.eps}
%%      \end{center}
%      \caption{A swallowtail}
%      \label{fig:swallowtail}
%   \end{minipage}
%\end{figure}

%\begin{figure}[!h]
%      \includegraphics[width=0.35\hsize]{cuspidal_edge.eps}
%      \caption{A cuspidal edge}
%      \label{fig:cuspidal_edge}
%\end{figure}

%
%
%\begin{center}
%\begin{figure}[!h]
%\begin{center}
%\begin{minipage}[t]{0.49\textwidth}
%\begin{center}
%         \includegraphics[width=0.35\hsize]{cuspidal_edge.eps}
%      \caption{A cuspidal edge}
%      \label{fig:cuspidal_edge}
%\end{center}
%\end{minipage}
%\hfill
%\begin{minipage}[t]{0.49\textwidth}
%\begin{center}
%         \includegraphics[width=0.6\hsize]{swallowtail.eps}
%      \caption{A swallowtail}
%      \label{fig:swallowtail}
%\end{center}
%\end{minipage}
%\end{center}
%\end{figure}
%\end{center}

We recall the criteria for cuspidal edges and swallowtails as in \cite{KRSUY}:

%\begin{df}
%A $C^\infty$ map $f$ on a 2-dimensional manifold $M^2$ to $\R^3$ is said the \textit{wave front}
%if the following propaties are satisfied:
%there exists a $C^\infty$ map $\nu:M^2 \longrightarrow \Sb^2$ ($\Sb^2$ is the unit sphere) such that
%(1) $\nu$ provides a unit normal vector field along $f$ except the set of the singular points and
%(2) $C^\infty$ map $(f,\nu):M^2 \longrightarrow \R^3 \times \Sb^2$ is an immersion.
%\end{df}

\begin{df}
Let $U \subset \R^2$ be an open domain.
A $C^\infty$ map $f: U \longrightarrow \R^3$ is called a \textit{frontal}
if there exists a $C^\infty$ map $\nu:U \longrightarrow \Sb^2$ ($\Sb^2$ is the unit sphere)
such that
$\nu(p)$ is perpendicular to $df(T_p U)$ for $p\in U$,
where $df$ is the differential of $f$ and
$T_p U$ is  the tangent space at $p$ to $U$.
Such a map $\nu$ is called a \textit{unit~normal~vector~field} of $f$.
Moreover, if the $C^\infty$ map $L:=(f,\nu):U \longrightarrow \R^3 \times \Sb^2$ is an immersion,
$f$ is called a \textit{wave front} (or \textit{front}).
\end{df}

%We prepare several definitions given in \cite{KRSUY}:
A singular point $p \in U$ of a frontal $f(u,v)$ is \textit{non-degenerate}
if the differential $d\lambda$ of $\lambda:=\det(f_u, f_v, \nu)$
%,which is $d\lambda=\lambda_u du + \lambda_vd v$,
does not vanish at $p$,
where $f_u:=\partial f / \partial u$ and $f_v:=\partial f / \partial v$.
If a singular point $p$ is non-degenerate,
the singular set of $f$ is a regular curve near $p$ on $U$.
This regular curve $c(s)$ is called the \textit{singular curve} of $f$, and
a tangent vector to $c$ is called a \textit{singular direction} of $f$.
Moreover, a nonzero vector $\eta \in T_{c(s)} U$ satisfying $df(\eta)=0$ is called
a \textit{null direction} of $f$.
We can take such $\eta(s)$ as a $C^\infty$ vector field along $c(s)$ near $p$, and
the $C^\infty$ vector field $\eta(s)$ is called a \textit{null vector field} of $f$.

%The following discriminant methods of types of singularities is known:
\begin{fact}[{\cite{KRSUY}}] \label{fact:cusp_swallow_dic}
Let $f=f(u,v):U \longrightarrow \R^3$ be a wave front.
We denote by $c(s)$
the singular curve near a non-degenerate singular point $p\in U$ of $f$
such that $c(0)=p$, and by $\eta(s)$ a null vector field along $c(s)$.
We set $\rho := \det(c', \eta)$. Then
\begin{enumerate}[$(i)$]
 \item $p$ is a cuspidal edge if and only if $\rho(0) \neq 0$,
 \item $p$ is a swallowtail if and only if $\rho(0)=0$ and $\rho'(0) \neq 0$.
\end{enumerate}
\end{fact}

Cuspidal edges and swallowtails are 
wave fronts as $C^\infty$ map germs.
It should be remarked that
criteria for cuspidal edges and swallowtails 
of developable surfaces have been given in \cite[Theorem 3.7]{Izumiya2}.
One can apply the criteria instead of those 
in Fact \ref{fact:cusp_swallow_dic}.

Next, we consider the a-completion of a ruled M\"obius strip $F(s,u)=\gamma(s)+u\xi(s)$ with singularities.
By a suitable change of parameters,
we may assume that $s$ is an arc-length parameter of $\gamma$ and $\xi(s)$ is a unit vector for $s\in \R$.

\begin{lem} \label{lem:Fs_cross_Fu}
Let $F(s,u)=\gamma(s)+u\xi(s)$ be a ruled M\"obius strip
where $s$ is arc-length and $\xi(s)$ is a unit vector for $s \in \R$.
%where $|\gamma'(s)|=|\xi(s)|=1$ for each $s$.
%Let $u$ move in $\R$.
Then
   \begin{equation*}
      \left| F_s \times F_u \right|^2
         =\left\{
            \begin{array}{ll}
               \left| \xi'(s) \right|^2
                 \left( u + \frac{\gamma'(s) \cdot \xi'(s)}{\left| \xi'(s) \right|^2} \right)^2
                 + \frac{\det(\gamma'(s), \xi(s), \xi'(s))^2}{\left| \xi'(s) \right|^2}
                  &\quad (\xi'(s)\neq \boldsymbol{0}) ,\\
             \left| \gamma'(s) \times \xi(s) \right|^2 &\quad (\xi'(s)= \boldsymbol{0}),
            \end{array}
         \right.
   \end{equation*}
   where the dot `$\,\cdot$' is the inner product and the cross `$\,\times\!$' is the vector product in $\R^3$.
\end{lem}

\begin{proof}
Since $F_s=\gamma' + u \xi'$ and $F_u=\xi$,
this assertion is obvious when $\xi'(s)=\boldsymbol{0}$. 
So we assume $\xi'(s)\neq \boldsymbol{0}$, and then
\begin{align*}
   \left| F_s \times F_u \right|^2
     &= \left| \xi' \right|^2
         \left( u + \frac{\gamma' \cdot \xi'}{\left| \xi' \right|^2} \right)^2
         + \frac{\left| \gamma' \times \xi \right|^2 \left| \xi' \right|^2 - (\gamma' \cdot \xi')^2}
                {\left| \xi' \right|^2}.
\end{align*}
Since
\begin{align*}
   |(\gamma' \times \xi) \times \xi'|^2
   &= |\gamma' \times \xi|^2 |\xi'|^2 - \left( (\gamma' \times \xi) \cdot \xi' \right)^2 
   = |\gamma' \times \xi|^2 |\xi'|^2 - \det(\gamma', \xi, \xi')^2
\end{align*}
and
\begin{align*}
   |(\gamma' \times \xi) \times \xi'|^2
   &= | (\gamma' \cdot \xi')\xi - (\xi \cdot \xi')\gamma' |^2
   = (\gamma' \cdot \xi')^2,
\end{align*}
we get the conclusion.
\end{proof}

%This lemma gives the following assertion.
%\begin{cor}
%A point $(s_0,u_0)$ is a singular point of a $C^\infty$ ruled surface $F(s,u)$ satisfying the conditions of the lemma above
%if and only if either $(1)$ or $(2)$ is satisfied;
%\begin{enum}
%\item $\xi'(s_0) \neq \boldsymbol{0}$,
%           $\displaystyle u_0 = -\frac{\gamma'(s_0) \cdot \xi'(s_0)}{\left| \xi'(s_0) \right|^2}$,
%           $\det(\gamma'(s_0), \xi(s_0), \xi'(s_0))=0$,
%      \vspace{0.5ex}
%\item $\xi'(s_0) = \boldsymbol{0}$,
%           $\gamma'(s_0) \times \xi(s_0) =\boldsymbol{0}$.
%\end{enum}
%\end{cor}

The equation $\det(\gamma', \xi, \xi')=0$ is
a necessary and sufficient condition of flatness of ruled M\"obius strips.
Hence, if $F$ is developable, the $C^\infty$ map
\begin{equation} \label{eq:nu}
   \nu(s,u) := \frac{\gamma'(s) \times \xi(s)}{|\gamma'(s) \times \xi(s)|}
\end{equation}
%$
%   \nu(s) := {\gamma'(s) \times \xi(s)}/{|\gamma'(s) \times \xi(s)|}
%$
gives a unit normal vector field along $F(s,u)$, so $F$ is a frontal.
Since $\nu$ does not depend on $u$ when $F$ is developable,
we regard as $\nu(s)=\nu(s,u)$ and denote $\nu'=\nu_s$.

%Such a frontal $F(s,u)$ is a front on an open neighborhood 
%of a singular point $(s_0,u_0)$ of $F$ if and only if $\nu'(s_0)\neq \boldsymbol{0}$.
%Lemma \ref{lem:Fs_cross_Fu} implies the following assertion:

Let `${\sim}$' be the equivalence relation which regards two points $(s,u)$ and $(s+l,-u)$ as the same point in $\R^2$,
where $l$ is the period of the closed curve $\gamma(s)$.
We set $M:=\R^2 / {\sim}$.
Then, $F$ can be regarded as a  $C^\infty$ map of $M$ into $\R^3$.

\begin{lem} \label{lem:singular_curve_graph_and_null_direction}
Suppose that $F(s,u)=\gamma(s)+u\xi(s)$ is a developable M\"obius strip, then
\begin{enumerate}[$(i)$]
   \item \label{enum:S(F)_non-d} each singular point of $F$ is non-degenerate,
   \item \label{enum:S(F)} the singular set $S(F)$ of $F$ is given by
      \begin{equation*}\label{eq:SingGraph1}
         S(F) = \left\{ (s,u) \in M \,;\, u=- \frac{|\gamma'(s) \times \xi(s)|^2}{\gamma'(s)\cdot \xi'(s)} ,\ \xi'(s)\neq \boldsymbol{0} \right\},
      \end{equation*}
   \item \label{enum:null_d} the null vector field of $F$ is given by
   \begin{equation*}\label{eq:NullDirection1}
%      \binom{1}{-\gamma'(s)\cdot \xi(s)}.
       \frac{\partial}{\partial s} -(\gamma'\cdot \xi) \frac{\partial}{\partial u}.
   \end{equation*}
\end{enumerate}
\end{lem}

\begin{proof}
%We remark that $\xi'(s)=\boldsymbol{0}$ if and only if $\gamma'(s) \cdot \xi'(s)=0$.
We set $\lambda:=\det (F_s, F_u, \nu)$.
Since the singular points of $F$ do not appear on asymptotic lines $u\mapsto (s_1,u)$ for $s_1 \in \R$ satisfying $\xi'(s_1)=\boldsymbol{0}$,
we have $\lambda_u = \gamma' \cdot \xi' / |\gamma' \times \xi| \neq 0$ on $S(F)$.
Therefore, we obtain $(\ref{enum:S(F)_non-d})$.
Since 
$
   \left| \gamma' \times \xi \right| \left| \xi' \right| = \left|\gamma' \cdot \xi'\right|
$
by flatness of $F$ and Lemma \ref{lem:Fs_cross_Fu},
we obtain $(\ref{enum:S(F)})$.
Let $(s_0,u_0)$ be a singular point.
Since $k:=\gamma'(s_0) \cdot \xi(s_0)\ (\neq 0)$ satisfies $F_s(s_0,u_0)=k F_u(s_0,u_0)$,
we have $(\ref{enum:null_d})$.
%$
%dF\left({\partial}/{\partial s} -k\, {\partial}/{\partial u} \right )=0
%$
%at $(s_0,u_0)$.
%$$
%dF\left(\frac{\partial}{\partial s} -k \frac{\partial}{\partial u} \right )=0.
%$$.
\end{proof}

Since $\xi$ is not a constant vector field,
there exists a point $s \in \R$ such that $\xi'(s) \neq \boldsymbol{0}$.
Therefore, the singular set $S(F)$ is not empty.
The following lemma gives a proof of the proposition in the introduction.

\begin{lem} \label{lem:exist_except_cusp}
Let $F(s,u)=\gamma(s)+ u\xi(s)$ be a developable M\"obius strip.
The a-completion of $F$ has at least one singular point other than cuspidal edge singularities
on each connected component of $S(F)$.
In particular, the a-completion of $F$ has at least one singular point other than cuspidal edge singularities.
\end{lem}

\begin{proof}
We remark that there exists a point $s \in \R$ such that $\xi'(s)=\boldsymbol{0}$,
since $\gamma' \cdot \xi'$ is an odd-periodic function.
Let $\{(s,u(s))\}_{s\in \R}$ be the graph of the singular curve of $F$ in the $(s,u)$-plane, and
let $\{(s,u(s))\}_{s_1<s<s_2}$ be a connected component of $S(F)$.
Then, the two points $s_1$ and $s_2$ satisfy $\xi'(s_1)=\xi'(s_2)=\boldsymbol{0}$ and 
$\xi'(s) \neq \boldsymbol{0}$ for $s \in (s_1,s_2)$.
Suppose $\gamma'(s)\cdot \xi'(s)>0$ for $s \in (s_1,s_2)$.
By Lemma \ref{lem:singular_curve_graph_and_null_direction} $(\ref{enum:S(F)})$, $u(s)$ satisfies
$$
   \lim_{s\searrow s_1} u(s) =\lim_{s\nearrow s_2} u(s) = -\infty,
$$
where $\searrow$ and $\nearrow$ mean approaching from above and below, respectively.
Then, the function
$$
   P(s):=-u(s)-\int_{s_1}^{s} \gamma'(t)\cdot \xi(t) dt
$$
satisfies
$$
   \lim_{s\searrow s_1} P(s) = \lim_{s\nearrow s_2} P(s) = \infty,
$$
since $|\gamma'(s)\cdot \xi(s)|<1$.
This implies that $P(s)$ attains a minimum at a point $s=s_0$.
Let $\rho(s)$ be the determinant of the $2 \times 2$ matrix consisting of
the two vectors for the singular direction and null direction of $F$.
Then, the function
$
      \rho(s) = -u'(s)-\gamma'(s)\cdot \xi(s) =P'(s)
$
vanishes at $s=s_0$.
By Fact \ref{fact:cusp_swallow_dic},
the singular point $(s_0,u(s_0))$ is not a cuspidal edge singularity.
The case $\gamma'\cdot \xi'<0$ is similar.
\end{proof}

\begin{rem} \label{rem:count_xi'_zoros}
Lemmas \ref{lem:singular_curve_graph_and_null_direction} and \ref{lem:exist_except_cusp} also imply that
the number of non-cuspidal-edge singular points on the a-completion of a developable M\"obius strip 
is greater than or equal to the number of connected components of the zero set of $\xi'$,
if these numbers are finite.

\end{rem}
%
%\begin{rem}
%%If $F$ is a front on a open domain $U$,
%the proof also give that $F$ must have cuspidal edges in $(s_1,s_2) \times \R$.
%Actually, if there do not exist any cuspidal edges, by Fact \ref{fact:cusp_swallow_dic}
%$\rho$ vanishes.
%This is a contradiction.
%\end{rem}

%\begin{proof}[Proof of Proposition]
%%Let us prove Propositon by using the lemmas above.
%Let $F(s,u)=\gamma(s)+u\xi(s)$ be a developable M\"obius strip.
%Since the odd-periodic function $\gamma'\cdot \xi'$ must have both of zeros and nonzeros,
%the open set
%$
%   \{s \in \R \;;\; \gamma'(s) \cdot \xi'(s) \neq \boldsymbol{0} \}
%$
%in $\R$ is not either the empty set or $\R$.
%We can choose a connected component of this open set and
%denote it by open interval $(s_1, s_2)$.
%Here, $s_1$ and $s_2$ may be $s_2=s_1+l$ if $\gamma'\cdot \xi$ has only one zero.
%By Lemma \ref{lem:exist_except_cusp}, we obtain the assertion of Propositon.
%\end{proof}

We close this section with an example having only one singular point other than cuspidal edge singularities.
This implies that the proposition gives the sharpest lower bound.

\begin{ex} We define a $2\pi$-periodic regular curve $\gamma = \gamma(s): \R \longrightarrow \R^3$ by
   \begin{equation*}
      \gamma(s):=\left(
                    \begin{array}{c}
                         \sin 2s   \\
                         \cos 2s   \\
                         ({1}/{\sqrt 2}) \sin s
                    \end{array}
                 \right),
%      \qquad (-\pi \leq t \leq \pi)
   \end{equation*}
   whose curvature function $\kappa(s)$ does not vanish.
   Let $\xi=\xi(s)$ be the $2\pi$-odd-periodic and non-vanishing vector field along $\gamma$ given by
   \begin{equation*}
      \xi(s):= p(s) \ve(s) + \cos ({s}/{2}) \vn(s) + \sin ({s}/{2}) \vb(s),
   \end{equation*}
   where $\ve$ is the unit tangent vector field, 
   $\vn$ is the normal vector field and
   $\vb$ is the binormal vector field of $\gamma$.
   Moreover, $\tau$ is the torsion function of $\gamma$ and
   $$
      p(s):= \frac{1}{\kappa(s)} \left.\left( \frac{1}{2|\gamma^\prime(s)|} + \tau(s)\right) \right/ \sin\frac{s}{2}.
   $$
%
%   Since
%   $$
%      \frac{1}{2|\gamma^\prime(0)|} = \frac{1}{3 \sqrt{2}} ,
%      \qquad \tau(0) = -\frac{1}{3 \sqrt{2}},
%   $$
   We remark that $p(s)$ and $\xi(s)$ are both smooth at $s=0, \pi$.
   Since $\det(\gamma', \xi, \xi') =0$, the map $F(s,u)=\gamma(s)+u\xi(s)$ is a developable M\"obius strip
   (See Figure~\ref{fig:1SingMobius_1}).

\begin{figure}[htbp]
 \begin{minipage}{0.49\hsize}
  \begin{center}
   \includegraphics[width=0.6\hsize]{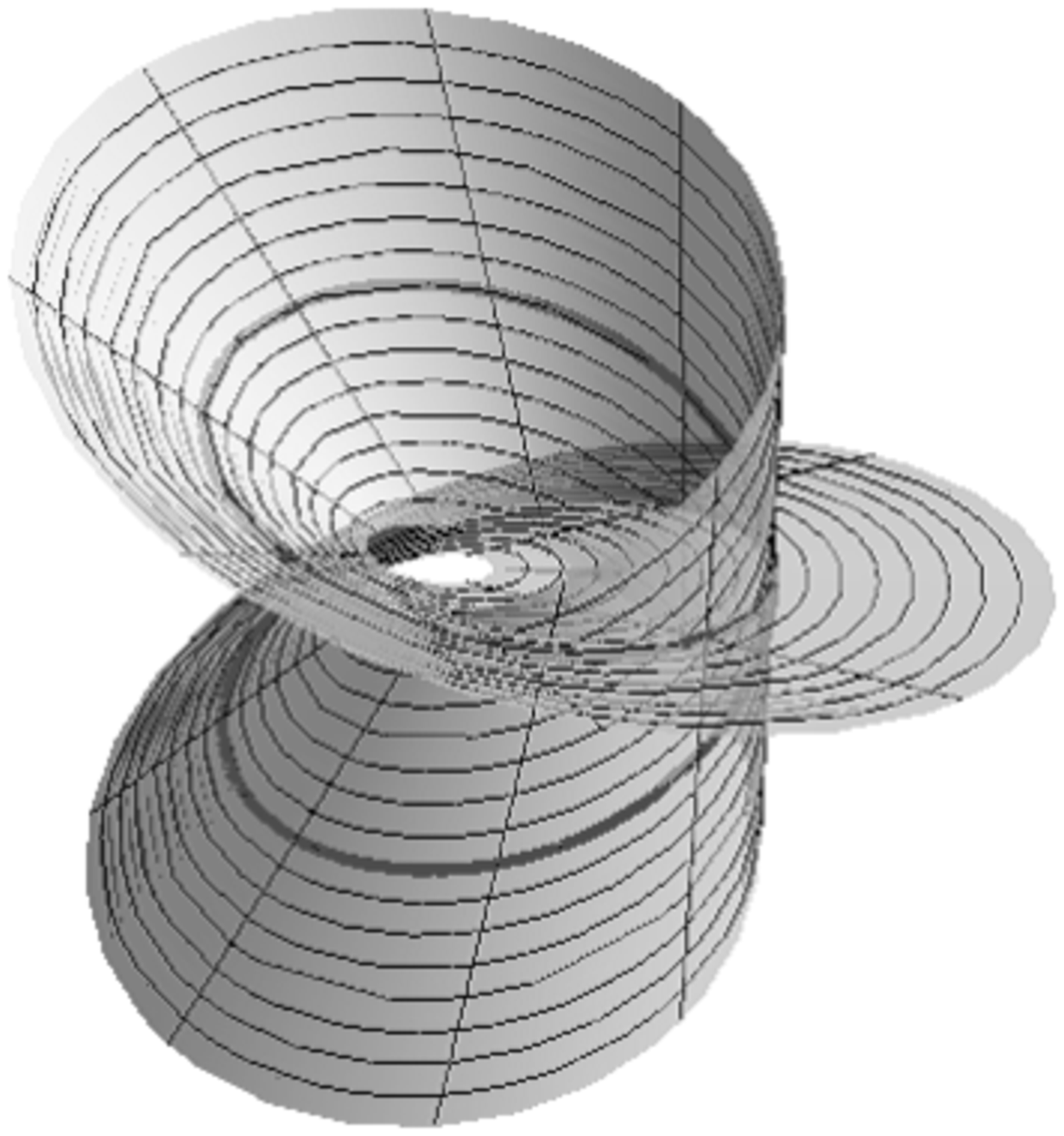}
  \end{center}
  \caption{{\small The image of $F(s,u)$}}
  \label{fig:1SingMobius_1}
 \end{minipage}
 \begin{minipage}{0.49\hsize}
  \begin{center}
   \includegraphics[width=0.55\hsize]{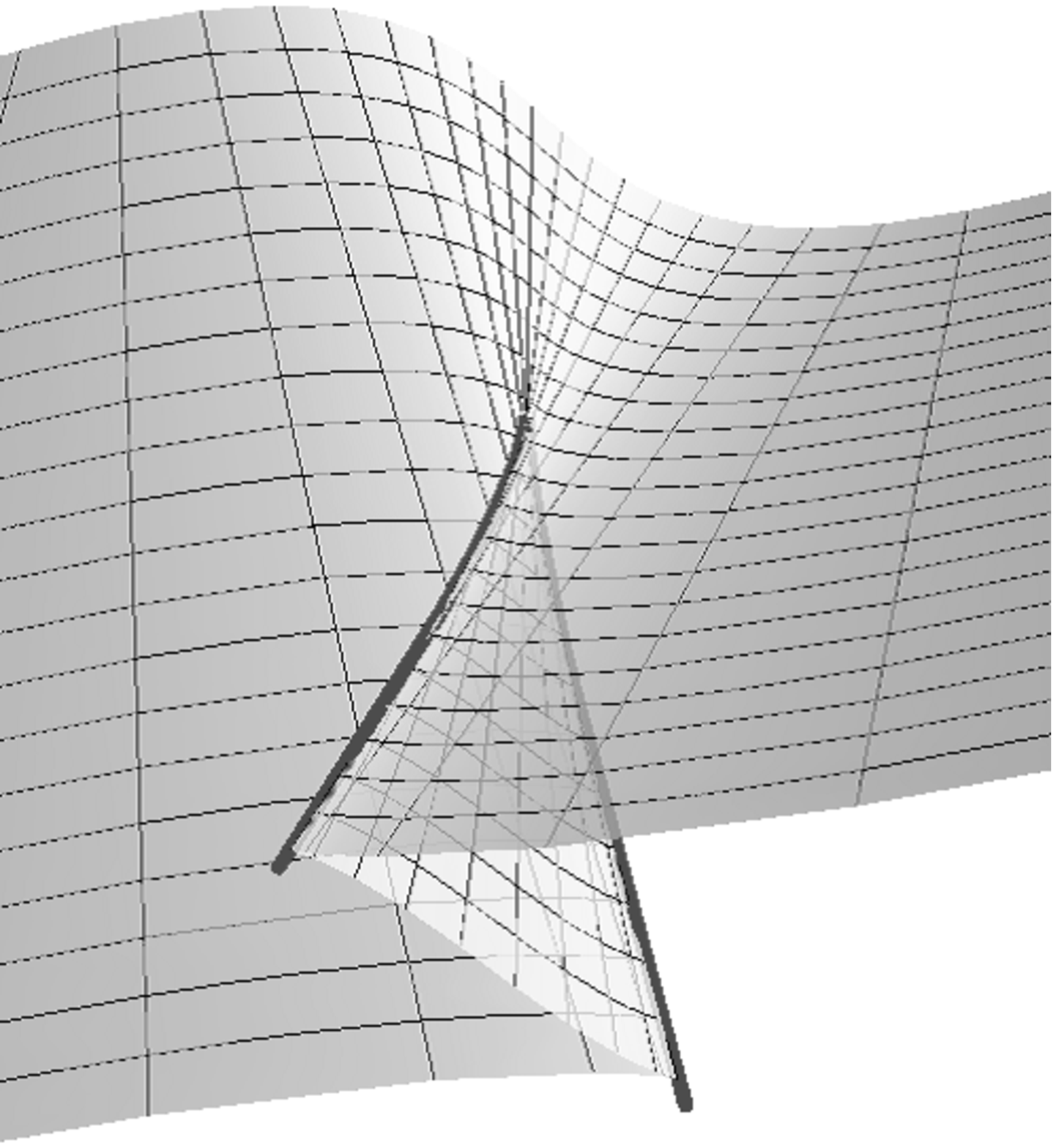}
  \end{center}
  \caption{{\small An open swallowtail}}
  \label{fig:open_swallowtail}
 \end{minipage}
\end{figure}

%
%\begin{figure}[htbp]
% \begin{minipage}{0.49\hsize}
%  \begin{center}
%   \includegraphics[width=0.8\hsize]{1SingMobius_1.eps}
%  \end{center}
%  \caption{}
%  \label{fig:1SingMobius_1}
% \end{minipage}
% \begin{minipage}{0.49\hsize}
%  \begin{center}
%   \includegraphics[width=0.7\hsize]{1SingMobius_2.eps}
%  \end{center}
%  \caption{}
%  \label{fig:1SingMobius_2}
% \end{minipage}
%\end{figure}
%
%
%
   The generating curve $\gamma(s)$ can be expressed by a rational function;
   if $x(s):=\tan (s/2)$ for $-\pi<s<\pi$, then we have
   $$
      \gamma(x)
               = \frac{1}{(1+x^2)^2}
                 \left(
                    \begin{array}{c}
                         4\; x (1-x^2)   \\
                         (1-2 x-x^2) (1+2 x-x^2)   \\
                         \sqrt{2}\; x(1+x^2)
                    \end{array}
                 \right).
   $$
   We set
   $$
      \hat\xi(x) := {\xi(s)}/{\cos (s/2)}
                  = \hat{p}(x)\ve(x) + \vn(x) + x \vb(x).
   $$
   Then $\hat{p}(x)$ is a $C^\infty$ function.
   Let $\rho(x)$ be the determinant of the $2 \times 2$ matrix
   consisting of the two vectors for the singular direction and null direction
   of $F(x,v)=\gamma(x)+v\hat\xi(x)$.
   We obtain\footnote{
   The software Mathematica (Version 7.0.0, Wolfram research) was used for this calculation.
   }
   \begin{equation*}\label{eq:lambda_1Sing}
      \rho = \frac{1}{|\hat\xi \times \hat\xi_x|^2}
                   \frac{\left(b_1 + b_{12} \sqrt{f_2}\right)^2}{\left(a_1  + a_{12} \sqrt{f_2}\right)}
                   \frac{x A(x)}
                        {%8 \sqrt{3}
                         (1+x^2)^3 (f_1)^{7/2} (f_2)^2 },
   \end{equation*}
   where $f_1(x) := 3 + 5 x^2 + 3 x^4$,\ $f_2(x) := 9 + 14 x^2 + 9 x^4$.
   Here, $a_1$, $a_{12}$, $b_1$, $b_{12}$ and $A$ are polynomials in $x$ 
   such that they have only even-degree terms and are non-negative.
   Moreover, the asymptotic line at $x=\infty$ has no singular points,
   so $\rho (x)=0$ if and only if $x=0$ (i.e.\ $s=0$).
   By Fact \ref{fact:cusp_swallow_dic}, 
   the singular point corresponding to $s=0$ is not a cuspidal edge singularity.
   On the other hand, $\nu'(s)=0$ if and only if $s=0$, where $\nu(s)$ is defined by \eqref{eq:nu}.
   Therefore, each singular point except at $s=0$ is a cuspidal edge, again by Fact \ref{fact:cusp_swallow_dic}.
\end{ex}

\begin{rem}

   By author's computer graphics,
   the singularity on the asymptotic line at $s=0$ looks like an `open swallowtail' (See Figure \ref{fig:open_swallowtail}; cf.\ \cite{Arnold}).
\end{rem}

\section{The proof of the theorem in the introduction}

Let $F(s,u)=\gamma(s)+u \xi(s)$ be a rectifying M\"obius strip (See the introduction).
Suppose that $s$ is an arc-length parameter of $\gamma$.
Then, $\gamma(s)$ and $\xi(s)$ satisfy
$$\gamma''(s) \cdot \xi(s)=0$$ for $s \in \R$,
since $\gamma(s)$ is a geodesic.
We normalize the ruling vector $\xi(s)$ for each $s \in \R$
such that
the projection of $\xi(s)$ into the rectifying plane at the point $\gamma(s)$
is a unit vector,
where the rectifying plane is a plane  perpendicular to $\ve(s)=\gamma'(s)$.
Then, $\xi$ can be expressed by
\begin{equation}\label{eq:normalized_Darboux_vector_field}
   D:=\frac{\tau}{\kappa} \ve + \vb
\end{equation}
when $\kappa$ is nonzero, where $\kappa$ is the curvature function, 
$\tau$ is the torsion function and $\{\ve, \vn, \vb\}$ is the Frenet frame of $\gamma$.
This vector field $D$ is called the \textit{normalized Darboux vector field} of $\gamma$.
The ratio $\sigma:=\tau/\kappa$ is called the \textit{conical curvature} of $\gamma$ (cf.\ Heil \cite{Heil}).

\begin{rem}
In \cite{Izumiya}, $D$ is called the {\it modified Darboux vector} 
along $\gamma$. Moreover, the criteria of cuspidal edges and 
swallowtails on the rectifying developable surfaces 
associated to $\gamma$ are given
in terms of conical curvature $\sigma=\tau/\kappa$. For example, 
$(s_0,u_0)$ is a non-cuspidal-edge singularity of
$F(s,u)=\gamma(s)+u D(s)$ if and only if
\red{$u_0 = -1/\sigma'(s_0)$}, $\sigma'(s_0)\neq 0$ and $\sigma''(s_0)=0$
(see \cite[Theorem 2.2]{Izumiya}). 
%One can also prove this criterion using 
%Fact \ref{fact:cusp_swallow_dic}.
\end{rem}

%and
%a extreme point of $\sigma$ is called a \textit{Darboux vertex}.}
%We remark that assuming that a closed curve satisfies \cite[Definition 1]{Heil} in $\R^3$, % (in particular, the curvature and torsion never vanish at the same points),
%he proved that the sum of three numbers of zeros of $\tau$, zeros of $\kappa$ and Darboux vertices is greater than or equal~to~four.
%
We recall the following facts in order to explain properties of the conical curvature
of a regular space curve.
%
%
%\begin{equation*}\label{eq:Heil}
%   V + K + D \geq 4,
%\end{equation*}
%where $V$ is the number of zeros of $\tau$, $K$ the number of zeros of $\kappa$ and $D$ the number of the Darboux vertices.
%We remark, however, that the assertion cannot apply to Theorem directly,
%since closed geodesics in rectifying M\"obius strips do not satisfy the definition.
%However,
%it was written that (\ref{eq:Heil}) was probably true under some assumptions, 
%even if $\kappa$ and $\tau$ vanish at the same points (See \cite[Remark 2]{Heil}).

\begin{fact}[{cf.\ \cite{Heil}}] \label{fact:sigma=kappa_g}
Let $I\subset \R$ be an open interval and
$\gamma:I \longrightarrow \R^3$ a regular curve.
If the curvature function $\kappa$ of $\gamma$ does not vanish,
then the geodesic curvature function of the unit tangent vector field $\ve:I \longrightarrow \Sb^2$ of $\gamma$ 
as a spherical curve is equal to the conical curvature $\sigma=\tau/\kappa$ of $\gamma$.
\end{fact}

A $C^\infty$ function $g=g(s):I \longrightarrow \R$ is said to be
\textit{strictly increasing} (resp.\ \textit{monotonically increasing in the wider sense})
if $g'(s) >0$ (resp.\ $g'(s) \geq 0$) for $s\in I$.
A regular spherical curve $\alpha=\alpha(s):I \longrightarrow \Sb^2$ is called
an \textit{honestly positive spiral} (resp.\ \textit{positive spiral})
if the geodesic curvature function of $\alpha$ is
strictly increasing (resp.\ monotonically increasing in the wider sense).

\begin{fact}[{\cite{Kneser,Pinkall}}] \label{fact:spiral_nest}
Let $\alpha=\alpha(s):I \longrightarrow \Sb^2$ be a honestly positive spiral $($resp.\ \textit{positive spiral}$)$.
We denote by $C(s) \subset \Sb^2$ the osculating circle of $\alpha$ at $s\in I$ and
assign $C(s)$ the orientation compatible with that of $\alpha(s)$
for each $s \in I$.
Let $D(s)$ be the left-hand domain of $C(s)$.
Then, $s_1<s_2$ implies
$\overline{D(s_2)} \subset D(s_1)$ $($resp.\ $D(s_2) \subset D(s_1)$$)$.
\end{fact}

We return to the initial settings in this section.
The Frenet frame of $\gamma(s)$ cannot be defined if $\kappa(s)= 0$.
However, we can construct an `extended' Frenet frame defined on $\R$ by using the ruling vector field $\xi(s)$.
We set
\begin{equation*}
      \hat\vn := -\frac{\ve \times \xi}{| \ve \times \xi |}, \quad 
      \hat\vb := \ve \times \hat\vn, \quad 
      \hat\kappa := \ve' \cdot \hat\vn, \quad 
      \hat\tau := -\hat\vb' \cdot \hat\vn, \quad
      \hat\sigma := \ve \cdot \xi.
\end{equation*}
These vector fields and functions are of class $C^\infty$.
Then $\{ \ve, \hat\vn, \hat\vb \}$ satisfies
$$
   \ve' = \hat\kappa \hat\vn, \qquad
   \hat\vn' = -\hat\kappa \ve + \hat\tau \hat\vb, \qquad
   \hat\vb' = -\hat\tau \hat\vn.
$$
Therefore, we have $\kappa = |\hat\kappa|$.
Moreover, if $\kappa(s) \neq 0$, then
$$
   \hat\vn(s)=\epsilon \vn(s), \quad \hat\vb(s)=\epsilon \vb(s),
$$
where $\epsilon := \hat\kappa(s)/\kappa(s)\ (=\pm 1)$.
The function $\hat\tau(s)$ is exactly equal to $\tau(s)$ if $\kappa(s) \neq 0$.
We can express the ruling vector field of $F$ by
$$
   \xi = \hat\sigma \ve + \hat\vb.
$$
Since $\det(\ve, \xi, \xi')=0$,
we have $\hat\tau=\hat\sigma \hat\kappa$.
Therefore, we regard $\hat\vn$, $\hat\vb$, $\hat\kappa$, $\hat\tau$ and $\hat\sigma$ as
smooth extensions of $\vn$, $\vb$, $\kappa$, $\tau$ and $\sigma$, respectively.
We set 
$$
    \hat K_+ := \{s \in \R \,;\, \hat\kappa(s)>0\}, \ 
    \hat K_0 := \{s \in \R \,;\, \hat\kappa(s)=0\}, \ 
    \hat K_- := \{s \in \R \,;\, \hat\kappa(s)<0\}.
$$
We regard $\ve=\gamma'$ as a closed curve in $\Sb^2$.
The spherical curve $\ve$ has singular points at zeros of $\kappa$.
For each $s\in\R$,
$\hat\vn(s)$ and $\hat\vb(s)$ can be regarded as
a unit tangent vector and a unit conormal vector of the spherical curve $\ve(s)$, respectively.
In particular, $\{ \hat\vn, \hat\vb, \ve \}$ gives a smooth positive orthonormal frame along $\ve$.
%
%
%
%Applying Lemma \ref{lem:singular_curve_graph_and_null_direction} for the rectifying M\"obius strip $F$,
%we obtain the following.
%
%\begin{lem} \label{lem:singular_curve_graph_and_null_direction_rectifying}
%The singular curve of the a-completion of
%a rectifying M\"obius strip $F(s,u)=\gamma(s)+u\xi(s)$
%is given in the $(s,u)$-plane by the graph 
%$
%   u = - 1/\hat\sigma'(s)
%$
%and the null direction of $F$ is given by
%$
%      \partial / \partial s.
%$
%\end{lem}
%
%\begin{proof}
%We write $\tilde\xi:=\xi/|\xi|$ and $v:=|\xi|u$.
%Then, we can apply Lemma \ref{lem:singular_curve_graph_and_null_direction} to $F=\gamma + v\tilde\xi$.
%\end{proof}
%
%Lemma \ref{lem:exist_except_cusp} and \ref{lem:singular_curve_graph_and_null_direction_rectifying}
%implies that
%the number of non-cuspidal-edge singular points on the a-completion of a rectifying M\"obius strip
%is greater than or equal to
%the number of the conected components of the zero point set of $\hat\sigma'(s)$
%if they are finite.
%
%
Since $\hat\sigma$ is of class $C^\infty$,
we can smoothly extend to $\R$ the osculating circle $C(s) \subset \Sb^2$ of $\ve(s)$.
In fact, the extended osculating circle $\hat C(s)$ can be canonically defined by a circle on $\Sb^2$
which passes $\ve(s)$ and whose center is 
$$
   \exp_{\ve(s)}\left(\frac{1}{2} \left(\arctan \frac{2}{\hat\sigma}\right) \hat\vb(s) \right),
$$
where $\exp_p:T_p \Sb^2 \longrightarrow \Sb^2$ is the exponential map at a point $p \in \Sb^2$.
%When $\hat\sigma(s)=0$, $\hat C(s)$ is a great circle on $\Sb^2$.
We assign $\hat C(s)$ the orientation compatible with the direction of $\hat\vn(s)$.
If $s\in \hat K_+$ (resp.\ $\hat K_-$), then the orientation of $\hat C(s)$ is equal (resp.\ opposite) to that of $C(s)$.
Let $\hat D(s)$ be the left-hand domain of $\hat C(s)$.

Since $\xi'=\hat\sigma' \ve$, by Remark \ref{rem:count_xi'_zoros},
it is sufficient to show that
the number of the connected components of the zero point set of $\hat\sigma'(s)$
is at least three.
We suppose that the number of locally maximal or locally minimal points of the odd-periodic function $\hat\sigma$ is only one.
We may assume that $s=0$ is the locally minimal point.
Then $\hat\sigma$ is a monotonic increasing function in the wider sense on the closed interval $[0,l]$, where $l$ is the period of $\gamma(s)$.
The restriction of the spherical curve $\ve$ to each connected component of $\hat K_+$ (resp.\ $\hat K_-$) is
a positive (resp.\ negative) spiral.
If we take two points $s_1$ and $s_2$ satisfying $s_1<s_2$ in each connected component of $\hat K_+ \cup \hat K_-$,
we have $D(s_2) \subset D(s_1)$ by Fact \ref{fact:spiral_nest}.
%Here, we do not take a closure to $D(s_2)$ as in Fact \ref{fact:spiral_nest}
%since $\hat\sigma$ may be constant on a interval.
%
On the other hand,
if we take two points $s_1$ and $s_2$ satisfying $s_1<s_2$ in each connected component of $\hat K_0$,
it holds that $\hat\kappa=\hat\tau=0$ on the closed interval $[s_1,s_2]$.
%   \begin{equation*}
%      \ve \cdot \hat\vn' = \ve' \cdot \hat\vn - (\gamma' \cdot \hat\vn)' =0, \quad
%      \hat\vn \cdot \hat\vn'= 0, \quad
%      \xi \cdot \hat\vn' = \xi' \cdot \hat\vn - (\xi \cdot \hat\vn)' =0,
%   \end{equation*}
%so $\hat\vn'=\boldsymbol{0}$.
Therefore $\hat\vn$ and $\hat\vb$ are constant on $[s_1,s_2]$, so
we have $D(s_2) \subset D(s_1)$.
Since the domain $D(s)$ depends smoothly on $s\in\R$,
we have $D(s_2) \subset D(s_1)$ for $s_1$ and $s_2$ satisfying $s_1<s_2$.
In particular, we obtain $D(l) \subset D(0)$.
On the other hand, 
the orientation of $C(l)$ is opposite to that of $C(0)$,
since $\hat\vn$ is odd-periodic.
Hence, $D(0) \cap D(l)$ is empty.
However, since $D(l)$ is not empty,
this is a contradiction.
Since $\hat\sigma$ is odd-periodic, 
$\hat\sigma$ must have
at least three locally minimal or locally maximal points.
Then, by Remark \ref{rem:count_xi'_zoros},
we obtain the theorem in the introduction.

\begin{ex}
We set
   \begin{equation*}
      \gamma(s):= \frac{1}{1 + (s+s^3)^2}
                  \left( \begin{array}{c}
                     ({2}/{5})s + s^3 +s^5 \\
                     s + s^3 \\
                     -{8}/{5} 
                  \end{array}\right),
   \end{equation*}
which gives a closed regular curve of $\Sb^1=\R \cup \{\infty\}$ in $\R^3$.
Moreover, $\gamma$ has only one inflection point at $s=\infty$.
We set $\hat\gamma(t):=\gamma(1/t)$.
Since
$$
   \hat\gamma'(t) \times \hat\gamma^{(3)}(t) |_{t=0} \neq \boldsymbol{0}, \quad
   \det(\hat\gamma'(t), \hat\gamma^{(3)}(t) , \hat\gamma^{(4)}(t) ) |_{t=0} =0
$$
and \cite[Corollary 2.11]{Kurono},
the $C^\infty$ map $F(s,u)=\gamma(s)+u\xi(s)$ is 
a rectifying M\"obius strip,
where $\xi(s)$ is as in \eqref{eq:normalized_Darboux_vector_field}.
The a-completion of $F$
has just three singular points other than cuspidal edges
(See Figures \ref{fig:3SingMobius} and \ref{fig:3SingMobius_1}).

\begin{figure}[htbp]
 \begin{minipage}{0.49\hsize}
   \begin{center}
      \includegraphics[width=0.7\hsize]{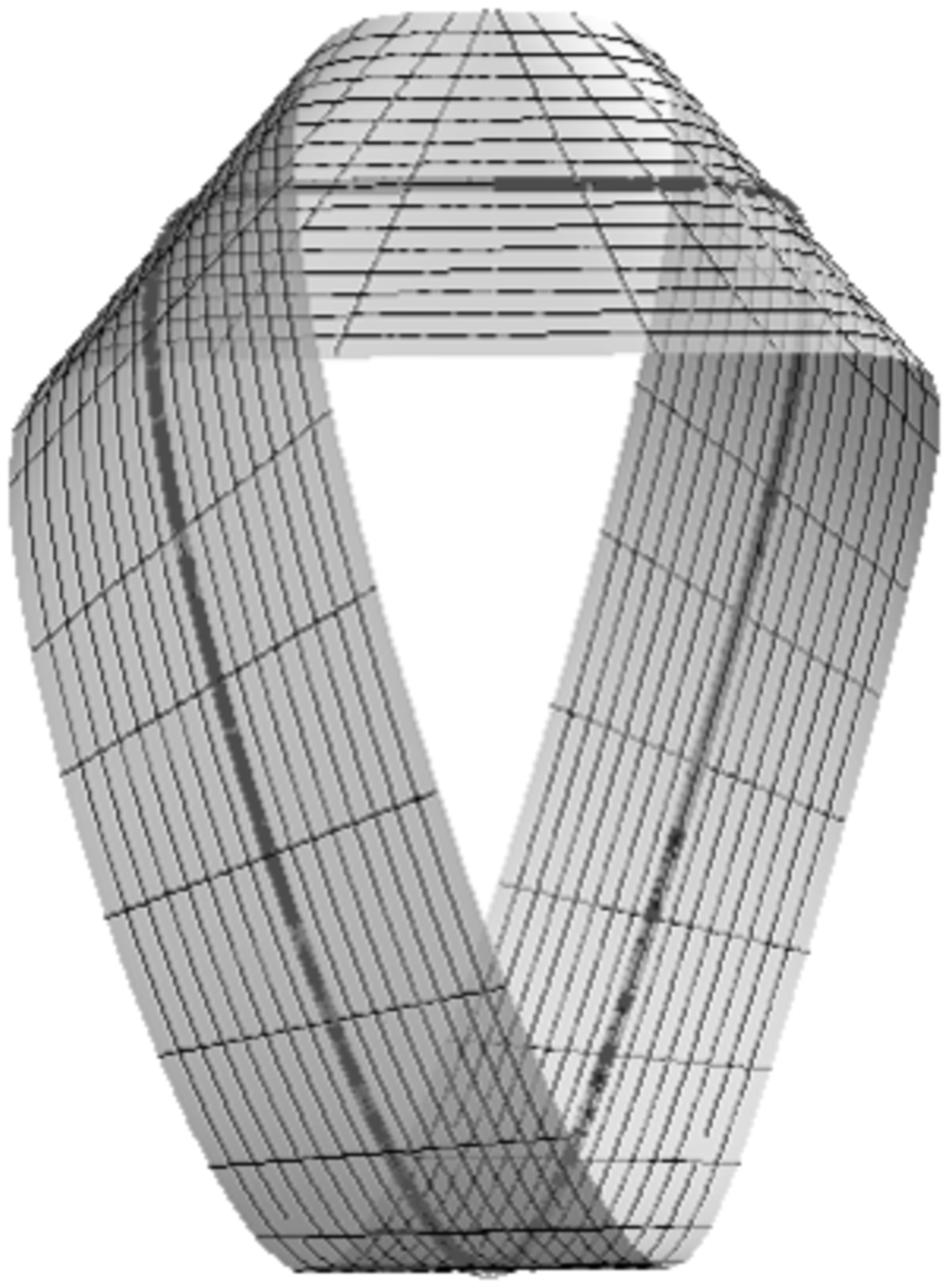}
   \end{center}
   \caption{{\small The image of $F(s,u)$}}
   \label{fig:3SingMobius}
 \end{minipage}
 \begin{minipage}{0.49\hsize}
  \begin{center}
   \includegraphics[width=0.7\hsize]{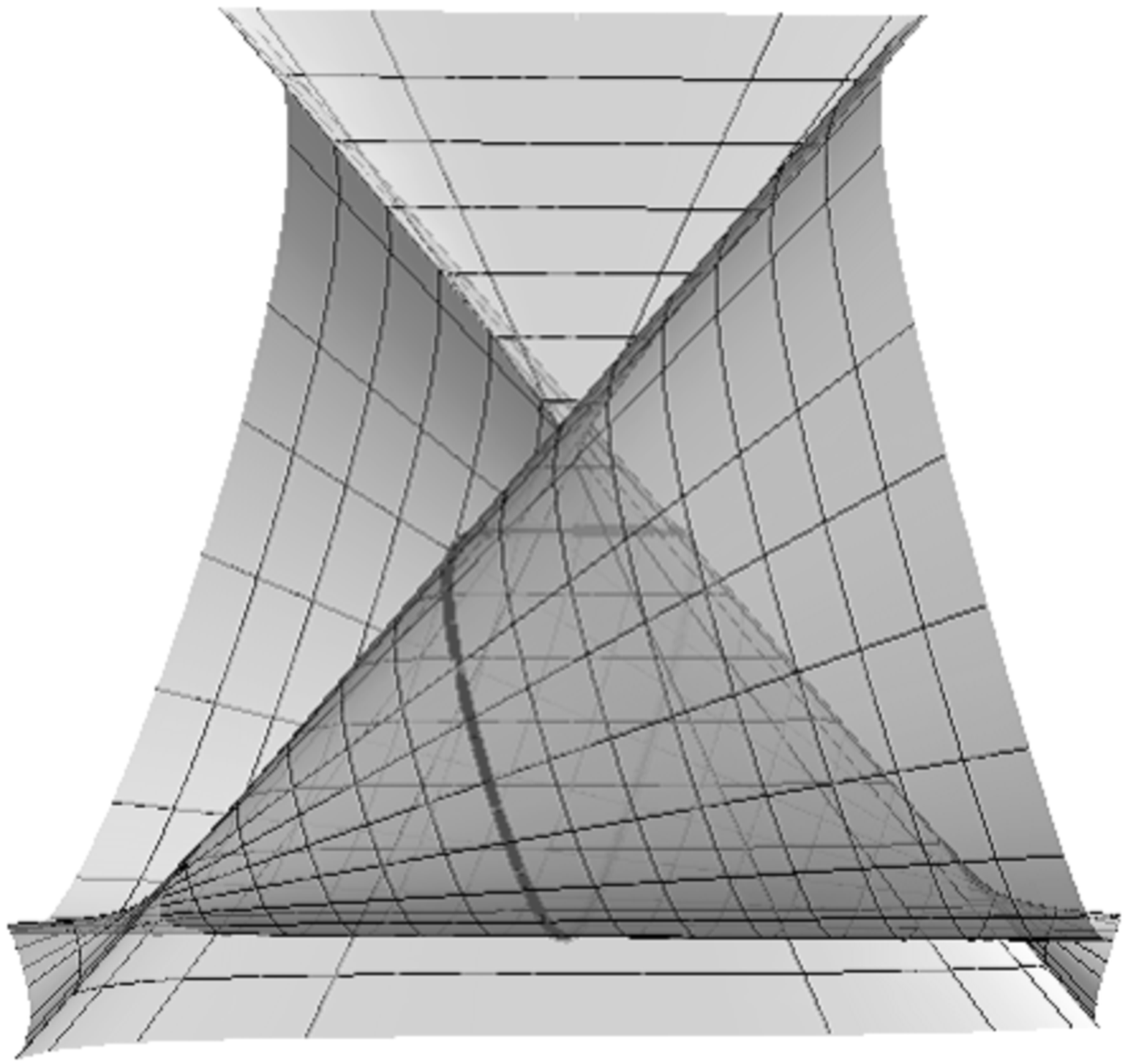}
  \end{center}
  \caption{{\small The~a-completion~of~$F$}}
  \label{fig:3SingMobius_1}
 \end{minipage}
\end{figure}

By Lemma \ref{lem:singular_curve_graph_and_null_direction},
the singular curve of $F$ is given by
$$
   s \longmapsto \left( s, u(s) := - \frac{|\gamma^\prime(s)|}{\hat\sigma'(s)}\right).
$$
Moreover, the null vector field of $F$ is ${\partial}/{\partial s}$.
We denote by $\rho(s)$
the determinant of the $2 \times 2$ matrix consisting of the two vectors for the singular direction and null direction of $F$.
Then, we can calculate\footnote{
   The software Mathematica (Version 7.0.0, Wolfram research) was used for this calculation.
   }
$$
      \rho(s):=u'(s)
      = \frac{(1 + s^2 + 2 s^4 + s^6) a(s)^{3/2} Q(s)}
             {s^2 b(s)^2},
$$
where $a(s)$, $b(s)$, $Q(s)$ are certain polynomials which have only even-dimensional terms and
$a(s), b(s)>0$.
It can be rigorously checked that
the polynomial $Q(s)$ has just two roots by Sturm's theorem.
Moreover, we have $\rho(s) \rightarrow 0\ (s\rightarrow \infty)$,
so $\rho(s)$ has three zeros including $s=\infty$.
On the other hand, since $\nu(s,u):=\hat\vn(s)$ is a unit normal vector field along $F$, 
the $C^\infty$ map $L=(F,\nu)$ is not immersed only at $(s,u)=(0,u(0))$.
Then,
$F$ has exactly three non-cuspidal-edge singular points by Fact \ref{fact:cusp_swallow_dic}.
We remark that the singularity at $(s,u)=(0,u(0))$ is a shape like an open swallowtail (See Figure \ref{fig:open_swallowtail}).
The other two non-cuspidal-edge singularities are both swallowtails.

\end{ex}

%\section*{Acknowledgements}
%\begin{acknowledgements}
%The author is grateful to Masaaki~Umehara for encouragement and suggestions.
%The author also thanks to Wayne~Rossman for carefully reading and commenting to this paper.
%Moreover, the author thanks to Yoshio Agaoka and Kentaro Saji for helpful information.
%\end{acknowledgements}

\medskip
\noindent
\textbf{Acknowledgements.}\ 
The author is grateful to Masaaki~Umehara for encouragement and suggestions.
The author also thanks to Wayne~Rossman for carefully reading and commenting to this paper.
Moreover, the author thanks to Yoshio Agaoka and Kentaro Saji for helpful information.


\medskip

\begin{thebibliography}{99}
   \bibitem{Arnold} V.~I.~Arnol'd: \textit{Lagrangian manifolds with singularities, asymptotic rays and the open swallowtail}, Funct.\ Anal.\ Appl. \textbf{15} (1981), 235--246.
   \bibitem{Chicone} C. Chicone and N. J. Kalton: \textit{Flat embeddings of the M\"obius strip in $\mathbb{R}^3$}, Commun.\ Appl.\ Nonlinear Anal. \textbf{9} (2002), 31--50.
   \bibitem{Heil} E. Heil: \textit{A four-vertex theorem for space curves}, Math.\ Pannon. \textbf{10} (1999), 123--132.
   \bibitem{Izumiya} S. Izumiya, H. Katsumi and T. Yamasaki: \textit{The rectifying developable and the spherical Darboux image of a space curve}, Banach Center Publications \textbf{50} (1999), 137--149.
   \bibitem{Izumiya2} S. Izumiya and N. Takeuchi: \textit{Geometry of ruled surfaces}, Applicable Mathematics in the Golden Age (edited by J. C. Misra), Narosa Publishing House, New Delhi, India (2003), 305--338
   \bibitem{Kneser} A. Kneser: \textit{Bermerkungen \"uber die Anzahl der Extreme der Kr\"ummung auf geschlossenen Kurven und \"uber verwandte Fragen in einer nichteuklidischen Geometrie}, Festschrift zum 70, Geburtstag von H.~Weber 1912, 170--180.
   \bibitem{KRSUY} M. Kokubu, W. Rossman, K. Saji, M. Umehara and K. Yamada: \textit{Singularities of flat fronts in hyperbolic 3-space}, Pacific J. Math. \textbf{221} (2005), 303--351.
   \bibitem{Kurono} Y. Kurono and M. Umehara: \textit{Flat M\"obius strips of given isotopy type in $R^3$ whose centerlines are geodesics or lines of curvature}, Geom.\ Dedicata \textbf{134} (2008), 109--130.
   \bibitem{Murata} S. Murata and M. Umehara: \textit{Flat surfaces with singularities in Euclidean 3-space}, J. Diff.\ Geom. \textbf{82} (2009), 279--316.
   \bibitem{Pinkall} U. Pinkall: \textit{On the four-vertex theorem}, Aequationes Math. \textbf{34} (1987), 221--230.
%   \bibitem{Randrup} T. Randrup, P. R{\o}gen, \textit{Sides of the Mobius strip}, Arch.\ Math. \textbf{66} (1996), 511--521.
%   \bibitem{Randrup} T. Randrup, P. R{\o}gen, \textit{How to twist a knot}, Arch.\ Math. \textbf{68} (1997), 252--264.
   \bibitem{Rogen} P. R{\o}gen: \textit{Embedding and knotting of flat compact surfaces in 3-space}, Comment.\ Math.\ Helv. \textbf{76} (2001), 589--606.
   \bibitem{Sabitov} I. K. Sabitov: Isometric immersions and embeddings of locally Euclidean metrics, Cambridge Scientific Publishers, 2009.
%   \bibitem{Sedykh} V. D. Sedykh, \textit{The four-vertex theorem of a convex space curve}, Funct.\ Anal.\ Appl. \textbf{26} (1992), 28--32.
   \bibitem{Wunderlich} W. Wunderlich: \textit{\"Uber ein abwickelbares M\"obiusband}, Monatsh.\ Math. \textbf{66} (1962), 276--289.
\end{thebibliography}
\end{document}